{\bf CURVES IN GRASSMANNIANS}
\bigskip
by Montserrat Teixidor i Bigas

\bigskip
{\bf INTRODUCTION}

Given a curve C and a line bundle L with sections, one can consider the
rational map
$$C \rightarrow {\bf P}(H^0(L))$$
 from the curve to projective space associated to the complete 
linear series of $L$.
A great deal is known about the image curve when the degree of the line
bundle is high enough. For example, the curve is 
projectively normal if $deg L \ge 2g+1$ (cf.[C],[M]).
 Saint Donat(cf[SD]) proved that the ideal of the curve is generated by
 quadrics and cubics and if $deg L\ge 2g+2$, it is generated by quadrics alone.
In the same vein, one can define condition $N_p$ for a curve. Roughly,
it means that the curve is projectively normal and in a minimal 
resolution of the ideal generators have the smallest possible degree
up to the $p^{-th}$ term.
Green (cf[G]) showed that if  $deg L \ge 2g+p+1$, then the curve satisfies 
property $N_p$.

Consider now a vector bundle $E$ of degree $d$ and rank $n\ge 2$ on the curve.
Assume that it has enough sections. We obtain a map
$$C\rightarrow G(n, H^0(E))$$
where G denotes the Grassmannian. In contrast with the line bundle
case, nothing seeems to be known about the
geometry of the image curve in such a Grassmanian, (say under good conditions
for $E$). In this paper, we try to clarify the situation somehow.

Our main results are in Propositions (1.1),(1.2). We consider the composition
$f$ of the following two maps
$$C\rightarrow G(n, H^0(E))\rightarrow {\bf P}(\wedge ^n H^0(E))$$
where the last map is the Plucker immersion of the Grassmannian in
projective space. In (1,1), we show that for very small degree the image of 
$f$ is a non-degenerate curve in projective space. In (1.2) we show that
for sufficiently high degree the image of $f$ lies on a proper linear
 subvariety and satisfies property $N_p$ as a curve in this variety.

The author is a member of the Europroj group ``Vector Bundles on Algebraic
Curves''. The author wants to thank the  Mathematics Department of Cambridge
 University where this paper was written. 
\bigskip

{\bf SET UP OF THE PROBLEM AND PROOF OF THE FIRST RESULTS}
\bigskip

We want to look at a curve from an external point of view. We need first 
decide on the set up.  Most interesting vector bundles are stable
 or at least semistable. In contrast with line bundles, even stable
 vector bundles of high degree don't behave all in the same way. We shall
restrict our attention to properties that hold for the generic vector 
bundle. As we want to define a map to the Grassmanian, we need
 $h^0(E)\ge n+2$. If the vector bundle is generic, this implies 
$d\ge ng+1$. 
 
The map $f$ is given by means of a 
(not necessarily complete) linear series. The line bundle that gives 
rise to this linear series is the determinant of the vector bundle E. 
The subspace of the space of sections is the image of the map 
$$\psi :\wedge ^nH^0(E)\rightarrow H^0(\wedge ^nE).$$
The first question that we want to ask is whether the curve is 
degenerate in projective space (i.e. lies on some hyperplane).
The answer is yes for sufficiently high degree. In fact the curve
is degenerate in projective space when the map $\psi$ above is 
not injective. For a generic vector bundle of degree
$d\ge n(g-1)$, $h^1(E)=0$. Therefore, $h^0(E)=d+n(1-g)$ and
 $dim\wedge ^nH^0(E)={d+n(1-g) \choose n}$. On the other hand
$\wedge ^n E$ is a line bundle and $h^0(\wedge ^n E)=d+1-g$.
Hence, when $d$ is large, $dim\wedge ^nH^0(E)\ge h^0(\wedge ^n E)$
and the map cannot be injective. For example, in rank two, injectivity
implies $d\le{4g-1+\surd {8g+1}\over 2}$.

We do not know whether injectivity holds in the whole allowable range
(for generic vector bundle). We prove in the next two paragraphs
that for very small degree, the map is indeed injective. On the other hand,
one can show that the map is surjective for generic $E$ if the degree is
high enough.
    The next two propositions  show our present state of knowledge.
One could ask in general whether $\psi$ always has maximal rank for a
generic $E$.
\bigskip

\proclaim (1,1) Proposition. Let E be a generic stable vector bundle
of degree ng+1 or ng+2. Then, the map $\psi$ is injective and therefore
the image by $f$does not lie on any proper linear subvariety.

\bigskip

\proclaim (1,2)Proposition. Let E be a generic stable vector bundle of degree 
$d\ge 2ng+1$. Then, the map $\psi$ is surjective. Hence the set of
hyperplanes in the ambient space of the Grassmanian cut a complete
linear system on the curve. If moreover $d\ge 2g+p+1$, then the image
curve satisfies property $N_p$ in the sublinear variety  that it spans.

Proposition (1,1) will be proved in the next two paragraphs.

 In order to
prove (1.2), it is enough to exhibit a vector bundle of degree d
 for which the map is surjective. We shall do this by induction on n.
 For rank one, there is nothing to prove. Assume the result true for
 vector bundles of rank $n-1$. Consider a generic extension
$$0\rightarrow L \rightarrow E\rightarrow F\rightarrow 0$$
where $deg L=2g$. Then $F$ is a generic vector bundle of rank 
$n-1$ and degree $deg F=deg E-2g\ge 2(n-1)g+1$.
Hence, the induction assumption applies to $F$ and the map
$$\wedge ^{n-1}H^0(F)\rightarrow H^0(\wedge ^{n-1}F)$$
is surjective.
We need to use the following Theorem of Mumford (cf.[M]).

\proclaim (1,3) Proposition (Mumford [M]).  If $L_1, L_2 $ are line bundles
$deg L_1\ge 2g, deg L_2 \ge 2g+1$, then the map
$$H^0(L_1)\otimes H^0(L_2)\rightarrow H^0(L_1\otimes L_2)$$
is surjective
\bigskip

By (1,3), the map 
$$H^0(L)\otimes H^0(\wedge ^{n-1}F)\rightarrow H^0(L\otimes \wedge ^{n-1}F)$$
 is surjective.
Moreover, as $degL=2g$, $h^1(L)=0$ and $H^0(E)=H^0(L)\oplus H^0(F)$.
Also $\wedge^n E=L\otimes \wedge ^{n-1}F$. The surjectivity of $\psi$
then follows.

The second statement follows directly from Green and Lazarsfeld`s Theorem.

\bigskip
{\bf (1,4) Remark} In characteristic zero, the hypothesis of genericity 
in (1,2) is unnecessary: It was proved by Butler (cf.[B]) that if 
$E, F$ are semistable vector bundles satisfying
$deg E=d>2g, \mu (F)\ge 2g$, then the map $H^0(E)\otimes H^0(F)\rightarrow
H^0(E\otimes F)$ is surjective. Hence, by induction on $n$, the map
$\otimes^n H^0(E)\rightarrow H^0(\wedge ^n E)$ is surjective. 
Moreover  the tensor product of semistable vector 
bundles is semistable. As $\otimes ^nE$ and $\wedge ^nE$ have the 
same slope, the symmetric product $S^n(E)$ is also semistable of
the same slope. Hence, $h^1(S^n(E))=0$ and the map $H^0(\otimes ^n(E))
\rightarrow H^0(\wedge ^n(E))$ is surjective
\bigskip

{\bf N+1 SECTIONS}

In this paragraph, we prove (1.1) for $d=ng+1$.
   Our method of proof is to degenerate the generic E to a direct sum 
of a generic line bundle $L$ of degree g and a generic vector bundle
 $F$ of rank $n-1$ and degree $(n-1)g+1$.
The map is not injective in this case but the elements in the kernel
cannot be deformed to a generic infinitessimal deformation of $E$. 
The techniques used are inspired of [W] Prop.(1.2).

With our assumptions,  $H^0(E)=H^0(L)\oplus H^0(F)$ and $h^0(L)=1$,
 $\wedge ^nH^0(E)=(H^0(L)\otimes \wedge^ {n-1} H^0(F))\oplus 
\wedge ^n H^0(F)$. 
 As $F$ is generic, it has $n$ sections and by induction on $n$,
 we can assume that the map

 $$\wedge  ^{n-1}H^0(F)\rightarrow H^0(\wedge ^{n-1}F).$$

is injective. Then, as $h^0(L)=1$, the map 

 $$ H^0(L)\otimes (\wedge ^{n-1}H^0(F))\rightarrow H^0(L\otimes \wedge ^{n-1}F)
=H^0(E)$$
is also injective. On the other hand, as $F$ has rank $n-1$
, $\wedge^n H^0(F)$ is contained in the kernel of $\psi$.
  Therefore, the kernel of $\psi$ in this case is 
$ \wedge ^n H^0(F)$.

Denote by $s^k,k=1...n$ a basis for $H^0(F)$. We want to see that 
$s^1\wedge ...\wedge s^n$ does not extend to an element in the kernel 
of $\psi _{\epsilon}$ for a generic infinitessimal deformation 
$E_{\epsilon}$ of $E$.
Recall that an infinitessimal deformation $E_{\epsilon}$ of $E$ is given by 
an element $e\in H^1(E^* \otimes E)$. Up to the choice of a suitable covering
 $C=\cup U_i$, we can represent $e$ by a cocicle
 $$\varphi _{ij} \in H^0(U_i\cap U_j,E^*\otimes E).$$
Then, $E_{\epsilon}$ can be represented locally as $E_{|U_i}\oplus \epsilon E_{\
|U_i}$
with gluings on $U_i\cap U_j$ given by the matrix
$$\pmatrix {Id&o\cr
\varphi_{ij}&Id}$$
In the case $E=L\oplus F$,
$$\varphi_{ij}=\pmatrix{ f^{11}_{ij}&f^{21}_{ij}\cr
                        f^{12}_{ij}&f^{22}_{ij}} $$
where $f^{kl}_{ij}\in H^0(U_i\cap U_j,Hom(L_k,L_l)),L_t=L(resp.F)$ if
 $t=1(resp.2)$.

Denote by
 $$s^k_{\epsilon}=(0,s^k)+\epsilon (s^k_{1i},s^k_{2i})$$

an extension to $E_{\epsilon}$ of the section $s^k$. Then
as $\epsilon ^2=0$ and $\wedge ^nH^0(F)\subset Ker \psi$, 

$$\psi_{\epsilon}(s^1_{\epsilon}\wedge...\wedge  s^n_{\epsilon})
=\epsilon \sum_{k=1}^n (-1)^ks^k_{1i}\wedge(\wedge _{l\not= k}s^l)$$
As F is a vector bundle of rank n-1, it is locally a direct sum of
n copies of the trivial bundle. When we take such a local representation,
the condition above can be written as

$$ det \pmatrix{s^1_{1i}&s^1\cr...&...\cr s^n_{1i}&s^n}=0$$
where each $s^i$ stands for a row of lenght n-1.
This means that the first column in the matrix is a linear combination of
 the remaining n-1. Equivalently, there exists a  local section on $U_i$
$\phi_i:F\rightarrow L$ such that $s^k_{1i}=\phi_i(s^k)$.
The gluing conditions can be read as $s^k_{2j}=s^k_{2i}+f^{22}_{ij}s^k$,
$s^k_{1j}=s^k_{1i}+f^{12}_{ij}s^k$. Replacing $s^k_{1i}=\phi _i (s^k)$ gives
$\phi_i+f^{21}_{ij}=\phi_j$ .

We now consider the set
$$\{ (\varphi _{ij},\phi_i,s^1_{2i}...s^n_{2i}) \in C^1(E^*\otimes E)\oplus
C^0(F^*\otimes L)\oplus( C^0 (F))^n, \phi_j=\phi_i+f^{21}_{ij},
s^k_{2j}=s^k_{2i}+f^{22}_{ij}s^k\}.$$

We can realize this set as the first hypercohomology group ${\bf H}$
 of the double complex of sheaves below:
$$\matrix{ C^0(E^*\otimes E)&\rightarrow &C^1(E^*\otimes E)&\rightarrow&
C^2(E^*\otimes E)\cr
\downarrow& &\downarrow& &\downarrow \cr
C^0(L\otimes F^*)\oplus( C^0 (F))^n&\rightarrow
&C^1(L\otimes F^*)\oplus C^1(F))^n&\rightarrow&
C^2(L\otimes F^*)\oplus( C^2 (F))^n}$$

The vertical maps in this diagram are induced by
 $(\pi^{21}, \pi^{22}(s^k))$.

We then have an exact sequence
$$0\rightarrow H^0(E^*\otimes E)\rightarrow
 H^0(L\otimes F^*)\oplus( H^0 (F))^n \rightarrow {\bf H}
\rightarrow H^1(E^*\otimes E)\rightarrow
 H^1(L\otimes F^*)\oplus( H^1 (F))^n$$
Note how this last map is defined by
$$\matrix {\alpha :&H^1(E^*\otimes E) &\rightarrow
& H^1(L\otimes F^*)\oplus( H^1 (F))^n\cr
&(\varphi _{ij})&\rightarrow&(f^{21}_{ij},f^{22}_{ij}(s^k))}$$

As $deg F=(n-1)g+1, h^1(F)=0$.
 As $deg L=g, h^1(L\otimes F^*)=(n-1)(g-1)+1\not= 0$,
$Ker \alpha =\{ \varphi _{ij}|f^{21}_{ij}=0\}$. Therefore, if the deformation
of $E$ does not preserve $F$ as a subundle,
  $s^1_{\epsilon}\wedge ... \wedge s^n_{\epsilon}\notin Ker \psi_{\epsilon}$
This completes the proof of the statement.

\bigskip
\bigskip
{\bf N+2 SECTIONS}

We now prove the analogous result for a generic vector bundle of
rank n with n+2 sections.

We first degenerate the vector bundle E to a direct sum  
$E=L^1_1\oplus F_2\oplus L^4_1\oplus ...\oplus L^n_1$ where subindices
 denote ranks
$deg L^1_1=g+1,degF_2=2g+1, deg L^4_1=...=degL^n_1=g$. We next degenerate
$F_2=L^2_1\oplus L^3_1, degL^2_1=g+1$. We then obtain
$$E=L^1_1 \oplus L^2_1 \oplus ... \oplus L^n_1, degL^1_1=degL^2_1=g+1, 
degL^3_1=...=degL^n_1=g$$
We assume all vector bundles appearing in the decompositions
above to be generic. Hence,
 $h^0(L^1_1)=h^0(L^2_1)=2,h^0(L^3_1)=...=h^0(L^n_1)=1$.
Then
$$\wedge ^nH^0(E)=\wedge ^2H^0(L^1_1)\otimes \wedge^2 H^0(L^2_1)\otimes 
(\oplus _{3\leq i_5<..< i_n}H^0(L^{i_5}_1)\otimes...
\otimes H^0(L^{i_n}_1))$$
$$\bigoplus \wedge^2 H^0(L^1_1)\otimes(\oplus _{i_k\neq 1}H^0(L^{i_3}_1)\otimes...
\otimes H^0(L^{i_n}_1))$$
$$\bigoplus \wedge^2 H^0(L^2_1)\otimes(\oplus _{i_k\neq 2}H^0(L^{i_3}_1)\otimes...
\otimes H^0(L^{i_n}_1))$$
$$\bigoplus H^0(L^1_1)\otimes...\otimes H^0 (L^n_1)$$ 
As $L_1^1$ is generic and $h^0(L^1_1)=2$, the kernel of the map
$H^0(L^1_1)\otimes H^0(L^2_1) \rightarrow H^0(L^1_1\otimes  L^2_1)$
is $H^0(L^2_1\otimes (L^1_1)^{-1})=0$. Moreover, as $H^0(L^i_1)=1,i=3..n$,
the composition map

 $$ H^0(L^1_1)\otimes...\otimes H^0 (L^n_1)\rightarrow 
H^0(L^1_1\otimes L^2_1)\otimes H^0(L^3_1)...\otimes H^0(L^n_1)
 \rightarrow H^0(L^1_1\otimes ...\otimes L^n_1)$$

is injective. So, the kernel of the map $\psi_{\epsilon}$ is given in this
 case by the other three summands in $\wedge ^nH^0(E)$.

We check first that the elements in 

$$\wedge ^2H^0(L^1_1)\otimes \wedge^2 H^0(L^2_1)\otimes
(\oplus _{3\leq i_5\leq..\leq i_n}H^0(L^{i_5}_1)\otimes...
\otimes H^0(L^{i_n}_1))$$

can be deformed to elements in the kernel of $\psi_{\epsilon}$ for any
infinitessimal deformation $E_{\epsilon}$ of $E$.

We choose basis
 $s^1,s'^1\in H^0(L^1_1),s^2,s'^2\in H^0(L^2_1), s^i\in H^0(L^i_1).$ 
We consider deformations of these sections to a generic $E_{\epsilon}$.
With the local description $E_{\epsilon}=E\oplus \epsilon E$ and 
taking into account that $E=L^1_1 \oplus L^2_1 \oplus ... \oplus L^n_1$,
the sections can be written as

$$s^1_{\epsilon}=s^1+\epsilon\overline s^1_i=(s^1,0,..,0)+
\epsilon (s^1_{1i},...,s^1_{ni})$$
$$...$$
$$s^n_{\epsilon}=s^n+\epsilon\overline s^n_i=(0,..,0,s^n)+
\epsilon (s^n_{1i},...,s^n_{ni}) $$ 
Independently of the choice of the $\overline s_i$,  we get
$$\psi _{\epsilon} (s^1_{\epsilon}\wedge s'^1_{\epsilon}\wedge s^2_{\epsilon}
\wedge s'^2_{\epsilon}\wedge s^{i_5}_{\epsilon}\wedge ...
\wedge s^{i_n}_{\epsilon})=0$$
We want to check now that these are the only elements that deform in 
an arbitrary direction as elements of the kernel. Assume the opposite. 
Then we can find an element in

$$ (*) \wedge^2 H^0(L^2_1)\otimes(\oplus _{i_k\neq 2}H^0(L^{i_3}_1)\otimes...
\otimes H^0(L^{i_n}_1))$$
$$\bigoplus \wedge^2 H^0(L^1_1)\otimes(\oplus _{i_k\neq 1}H^0(L^{i_3}_1)\otimes...
\otimes H^0(L^{i_n}_1))$$

in the kernel of $\psi_{\epsilon}$ that deforms in an arbitrary direction.
We first consider an element of the form
$s^1\wedge s'^1 \wedge s^{i_3}\wedge...\wedge s^{i_n},\{ i_3,..,i_n\}\cup 
\{k\}=\{ 2,..,n\}$. 
We want to show that these elements deform in the infinitessimal directions
$\varphi _{ij}=(f^{pq}_{ij})$ such that  $f^{pq}_{ij}=0$ if $p=1,q=k$.
Note that
 $\psi_{\epsilon} (s^1\wedge s'^1 \wedge s^{i_3}\wedge...\wedge s^{i_n}) 
=\epsilon( s^1s'^1_{ki}-s^1_{ki}s'^1)s^{i_3}...s^{i_n}$. 
This expression will be zero if and only if $s^1s'^1_{ki}-s^1_{ki}s'^1=0$.
As $s^1,s'^1$ have no common zeroes,$s^1_{ki}=s^1t_i, s'^1_{ki}=s'^1t_i$. 
Define
 $$A= (L^1_1)^n\oplus...\oplus( L^k_1(L^1_1)^{-1})\oplus(L^k_1)^{n-2}\
\oplus ...\oplus  (L^n_1)^n$$
We now consider the set
$$\{ (\varphi _{ij},s^1_{1i},s'^1_{1i}..s^n_{1i},...,t_i,s^3_{ki}
..s^n_{ki}...s^n_{ni}) \in C^1(E^*\otimes E)\oplus
C^0(A), t_j=t_i+f^{1k}_{ij},s^l_{mj}=s^l_{mi}+f^{lm}_{ij}s^l\}.$$

We can realize this set as the first hypercohomology group ${\bf H}$
 of the double complex of sheaves below:
$$\matrix{ C^0(E^*\otimes E)&\rightarrow &C^1(E^*\otimes E)&\rightarrow&
C^2(E^*\otimes E)\cr
\downarrow& &\downarrow& &\downarrow \cr
C^0(A)&\rightarrow&C^1(A)&\rightarrow&C^2(A)}$$ 
We then have an exact sequence
$$0\rightarrow H^0(E\otimes E^*)\rightarrow H^0(A)\rightarrow {\bf H}
\rightarrow H^1(E\otimes E^*)\rightarrow H^1(A)$$
We are interested in the kernel of this last map. Note that $h^1(L^i_1)=0,
h^1(L^k_1\otimes (L^1_1)^{-1})=g$. Hence, the kernel consists of those 
$\varphi _{ij}$ such that $f^{1k}_{ij}=0$.

Assume now that an element in (*) deforms in a generic direction
Write this element as 
$$\sum a_{i_3...i_n}s^1\wedge s'^1\wedge s^{i_3}\wedge ...\wedge s^{i_n}+
\sum b_{i_3...i_n}s^2\wedge s'^2\wedge s^{i_3}\wedge ...\wedge s^{i_n}$$
 Then, it 
deforms also in the directions that have $$f^{1k}_{ij}=0,k=2,...,n,
f^{2k}_{ij}=0, k\not= 2,l, f^{2l}_{ij}\not= 0.$$
Then, the coefficient $b_{i_3,...,i_n}=0$. Reasoning in this way for all l 
and replacing then 2 by 1, we obtain that all coefficients are 0. In particular
this completes the proof for $n=2,3$.

We now go back to the degeneration
 $E=L^1_1\oplus F_2 \oplus L^4_1 \oplus...\oplus L^n_1$. 
We need only show that the elements in 

$$\oplus_{\{ i_5,...,i_n\} \subset \{ 4,...,n\}}\wedge^2 H^0(L^1_1) \otimes
\wedge ^2 H^0(F_2)\otimes H^0(L^{i_5}_1)\otimes ...\otimes H^0(L^{i_n}_1)$$
$$\bigoplus$$
$$\oplus_{\{ i_6,...,i_n\} \subset \{ 4,...,n\}}\wedge^2 H^0(L^1_1) \otimes
\wedge ^3 H^0(F_2)\otimes H^0(L^{i_6}_1)\otimes ...\otimes H^0(L^{i_n}_1)$$

are not in the kernel. 
Reasoning as before, we see that the the only sections that deform
to sections in the kernel are those in

 $$\oplus_{\{ i_6,...,i_n\} \subset \{ 4,...,n\}}\wedge^2 H^0(L^1_1) \otimes
\wedge ^3 H^0(F_2)\otimes H^0(L^{i_6}_1)\otimes ...\otimes H^0(L^{i_n}_1)$$ 
Therefore, the only sections that could originally be deformed in all
directions are those in 
$$\oplus_{\{ i_6,...,i_n\} \subset \{ 4,...,n\}}\wedge^2 H^0(L^1_1) \otimes
\wedge ^2 H^0(L^2_1) \otimes  H^0(L^3_1)\otimes  H^0(L^{i_6}_1)\otimes ...
\otimes H^0(L_1^{i_n})$$

By symmetry in the $L^i, i=3,...,n$, this concludes the proof. 

\bigskip
\bigskip
{\bf REFERENCES}
\bigskip
[B] D.Butler, Normal generation of vector bundles over a curve,J.Diff
.Geom.{\bf 39} (1994), 1-34.

[C] G.Castelnuovo. Sui multipli di una serie di gruppi di punti appartenente
ad una curva algebrica, Rend.Circ.Mat.Palermo {\bf 7}(1892), 99-119.

[G] M.Green. Koszul cohomology and the geometry of projective varieties.
J.Diff.Geom. {\bf 19}(1984), 125-171.

[M] D.Mumford Varieties defined by quadratic equations. Corso CIME. In 
questions on algebraic varieties, Rome 1970. 30-100.

[S-D] B.Saint Donat. Sur les equations definissant une courbe algebrique.
C.R.Acad.Sc. Paris {\bf t 274} (1972), 324-327.

[W] G.Welters. Polarised abelian varieties and the heat equation. Comp.
Math. {\bf 49}(1983), 173-194.

\bigskip
Mathematics Department, Tufts University, Medford MA 02155, U.S.A.

temporary address DPMMS, 16 Mills Lane, Cambridge CCB2 1SB, England

teixidor@dpmms.cam.ac.uk

\end